\documentclass{iopart}

\usepackage{setstack}
\usepackage{amssymb}
\usepackage{graphicx}
\usepackage{amsthm}
\usepackage{a4wide}
\usepackage{harvard}

\eqnobysec
\def\beql{\smallskip
\begin{equation}}
\def\eeql{\smallskip
\end{equation}}
\def\beqal{\smallskip
\begin{eqnarray}}
\def\eeqal{\smallskip
\end{eqnarray}}

\def\beq*{\smallskip
\begin{equation*}}
\def\eeq*{\smallskip
\end{equation*}}
\def\beqa*{\beq*
\begin{array}{ll}}
\def\eeqa*{\end{array}
\eeq*}

\def\d{\textup{d}}
\def\pd{\partial}

\def\and{\quad\text{and}\quad}

\def\leqs{\leqslant}
\def\geqs{\geqslant}

\newtheorem{thm}{Theorem}[section]
\newtheorem{lem}[thm]{Lemma}
\newtheorem{prop}[thm]{Proposition}

\newtheorem*{defn}{Definition}

\begin{document}

\title{A Bernoulli linked-twist map in the plane}

\author{James Springham$^1$ and Stephen Wiggins$^2$}

\address{$^1$School of Mathematics, University of Leeds, Leeds LS2 9JT, United Kingdom\\$^2$School of Mathematics, University of Bristol, Bristol BS8 1TW, United Kingdom}

\eads{\mailto{j.springham@leeds.ac.uk}, \mailto{s.wiggins@bristol.ac.uk}}

\begin{abstract}
We prove that a Lebesgue measure-preserving linked-twist map defined in the plane is metrically isomorphic to a Bernoulli shift (and thus strongly mixing). This is the first such result for an explicitly defined linked-twist map on a manifold other than the two-torus. Our work builds on that of \citeasnoun{woj} who established an ergodic partition for this example using an invariant cone-field in the tangent space.
\end{abstract}

\address{Keywords: Linked-twist map, Mixing, Bernoulli property}


\section{Introduction}\label{Section1}
\setcounter{equation}{0}

\begin{figure}[htp]
\centering
\includegraphics[totalheight=0.25\textheight]{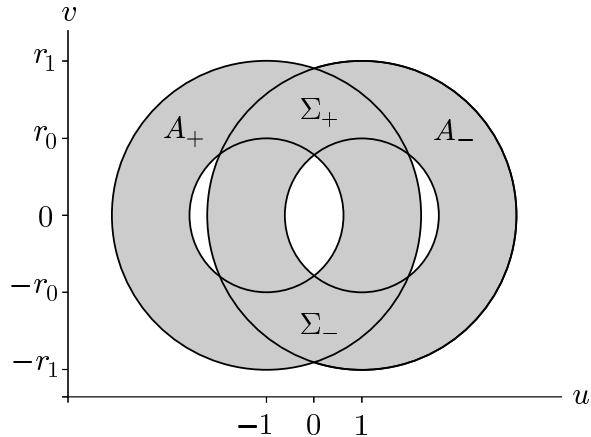}
\caption[The manifold $A\subset\mathbb{R}^2$]{The manifold $A\subset\mathbb{R}^2$ (shaded).}
\label{fig:A_plane}
\end{figure}

Let $\mathbb{S}^1=[-\pi,\pi]$ with opposite ends identified, fix $0<r_0<r_1<\pi$ and define an annulus
\beq*
L=\{(r,\theta):r_0\leqslant r\leqslant r_1,\theta\in\mathbb{S}^1\}.
\eeq*
The functions $M_{\pm}:\mathbb{R}^+_0\times\mathbb{S}^1\to\mathbb{R}^2$ given by
\beq*
M_{\pm}(r,\theta)=\pm(r\cos\theta-1,r\sin\theta)
\eeq*
map $L$ into the plane, the images $A_{\pm}=M_{\pm}(L)$ being centred at $(-1,0)$ and at $(1,0)$ respectively. Let $(u,v)$ denote the usual Cartesian coordinates in $\mathbb{R}^2$. We assume that $r_0,r_1$ are such that the annuli intersect in two disjoint regions, and denote the intersection region in which the $v$ coordinate is positive by $\Sigma_+$ and the other by $\Sigma_-$. See Figure~\ref{fig:A_plane}.

\begin{figure}[htp]
\centering
(a)\includegraphics[totalheight=0.13\textheight]{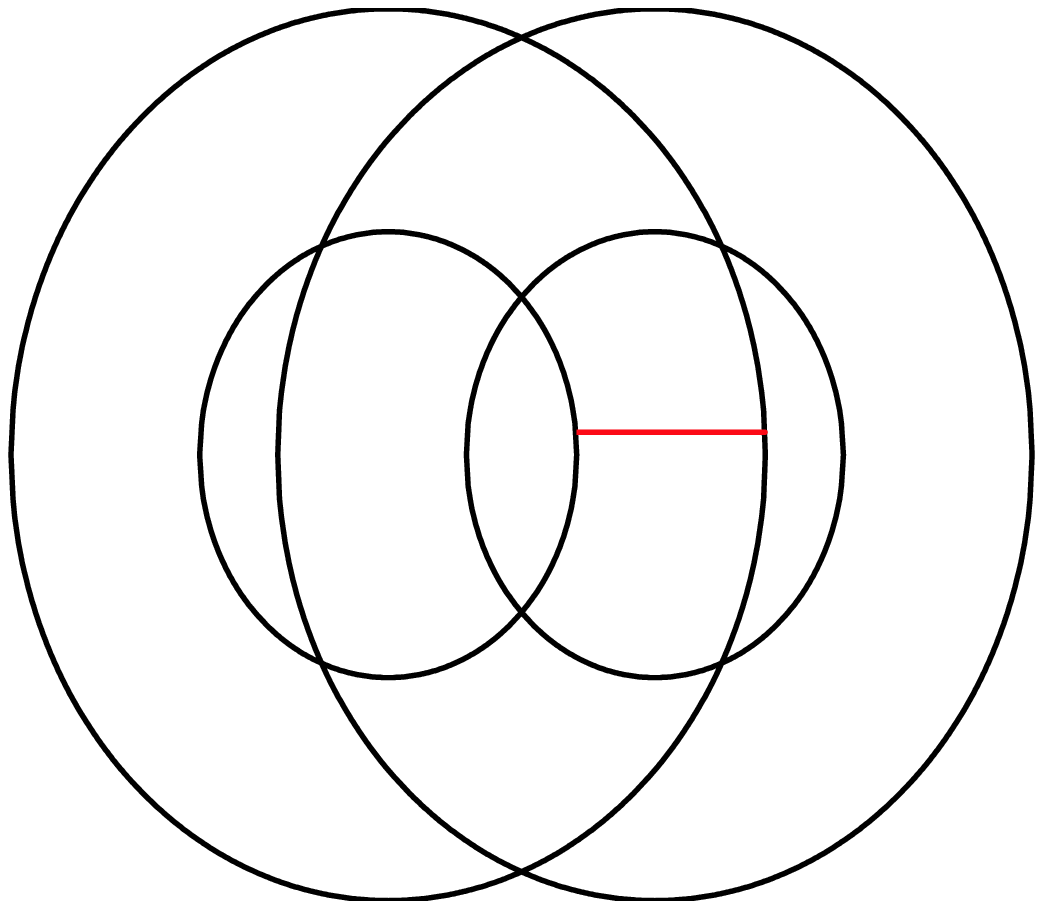}
(b)\includegraphics[totalheight=0.13\textheight]{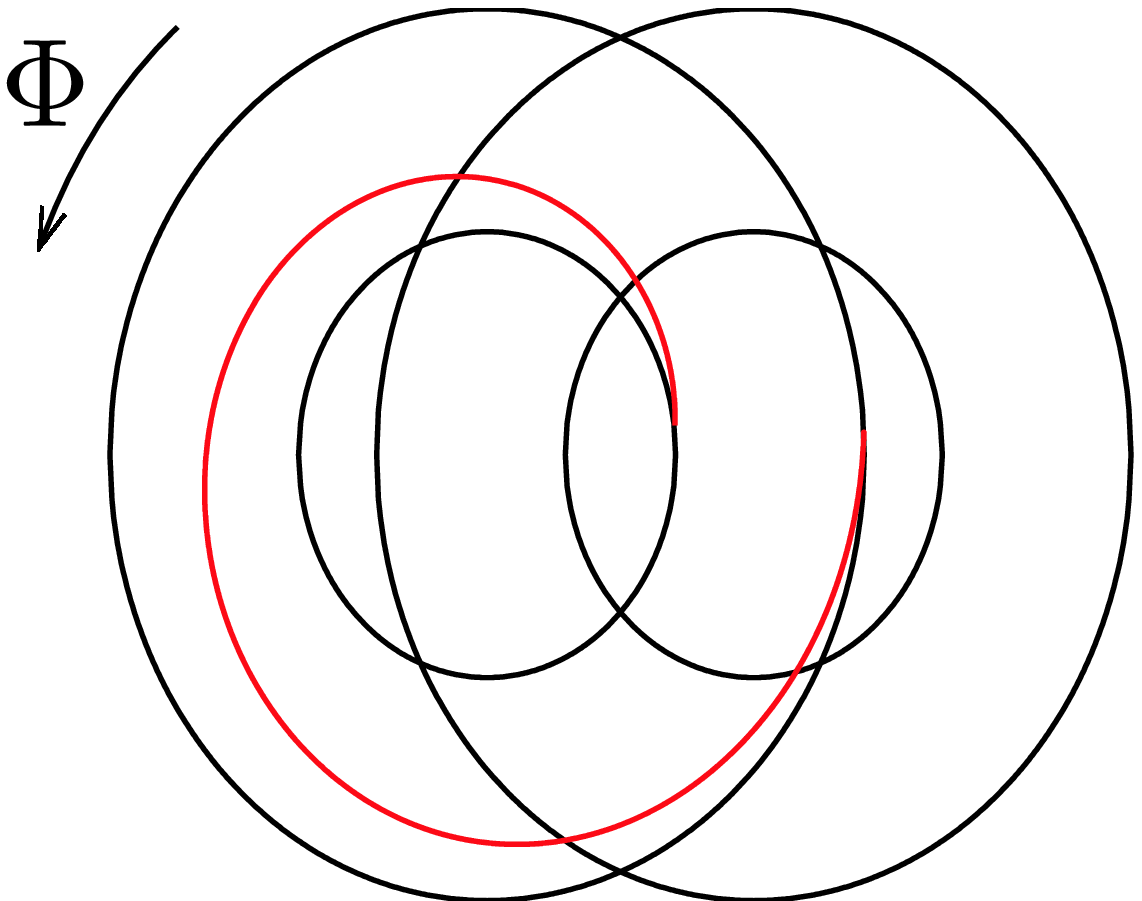}
(c)\includegraphics[totalheight=0.13\textheight]{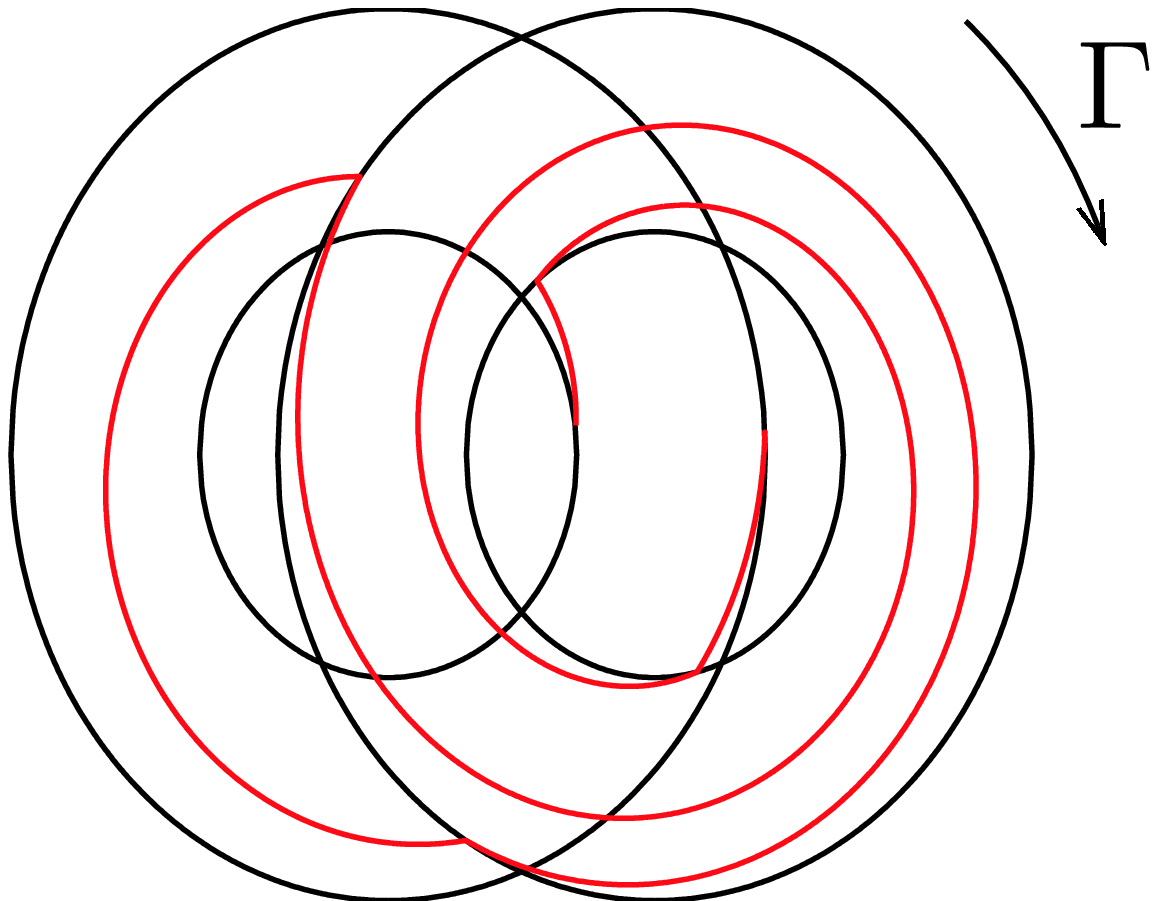}
\caption[A planar linked-twist map]{One iteration of the planar linked-twist map. Part~(a) shows some initial conditions in the form of a red horizontal line across the left-hand annulus $A_+$. Part~(b) shows the image of these points under the twist map $\Phi$ and part~(c) shows the image under the linked-twist map $\Theta=\Gamma\circ\Phi$.}
\label{fig:planar_twist}
\end{figure}

We denote $A=A_+\cup A_-$ and $\Sigma=\Sigma_+\cup\Sigma_-$. Inverses $M_{\pm}^{-1}:A\to\mathbb{R}_0^+\times\mathbb{S}^1$ to the above functions are given by
\beq*
M_{\pm}^{-1}=\left(\sqrt{(1\pm u)^2+v^2},\tan^{-1}\frac{v}{u\pm 1}\right).
\eeq*
Define a \emph{twist map} $\Lambda:L\to L$ by
\beq*
\Lambda(r,\theta)=(r,\theta+2\pi(r-r_0)/(r_1-r_0)),
\eeq*
remarking that $\Lambda$ leaves invariant the boundaries of $L$ and otherwise rotates points about the origin by an angle that increases with the radial coordinate. The \emph{twist function} $r\mapsto 2\pi(r-r_0)/(r_1-r_0)$ has derivative $c=2\pi/(r_1-r_0)$ and is affine. We define two twist maps $\Phi,\Gamma:A\to A$ in the plane, given by
\beql\label{eqn:phi}
\Phi(u,v)=\left\{
\begin{array}{r@{\quad}l}
M_+\circ\Lambda\circ M_+^{-1}(u,v) &\text{if }(u,v)\in A_+ \\
(u,v) &\text{otherwise,}
\end{array}
\right.
\eeql
\beql\label{eqn:gamma}
\Gamma(u,v)=\left\{
\begin{array}{r@{\quad}l}
M_-\circ\Lambda^{-1}\circ M_-^{-1}(u,v) &\text{if }(u,v)\in A_- \\
(u,v) &\text{otherwise.}
\end{array}
\right.
\eeql
\begin{defn}[Linked-twist map]
A linked-twist map $\Theta:A\to A$ is the composition $\Theta=\Gamma\circ\Phi$.
\end{defn}
We call this a linked-twist map \emph{in the plane} or say that it is \emph{planar}. It preserves the Lebesgue measure $\mu$ on $A$.  The relative direction of the two rotations (here they are opposite) affects the ergodic properties of linked-twist maps. \citeasnoun{sturman} discuss this in some detail; in their terminology the present map is \emph{co-twisting}. We illustrate $\Theta$ in Figure~\ref{fig:planar_twist}.

The purpose of this paper is to prove the following:
\begin{thm}
\label{thm:main_plane}
Let $r_0=2$ and $r_1=\sqrt{7}$. The planar linked-twist map $\Theta:A\to A$ is isomorphic to a Bernoulli shift.
\end{thm}

We make some remarks. Crucial progress toward this result was made by \citeasnoun{woj} and our work builds upon his. He proved that the system considered here is amongst a family of such systems that possess an \emph{ergodic partition}, using the technique of finding an invariant cone-field; we review this result in Section~\ref{Section3}. The work lead him to conjecture that such systems also are mixing and our result shows that this is indeed the case for the system considered. We discuss the reasons for restricting to this one example in Section~\ref{Section8}.

Our paper is organised as follows. We discuss the recent resurgence of interest in linked-twist maps in Section~\ref{Section2} and the work of Wojtkowski in Section~\ref{Section3}. The cornerstone of our proof of Theorem~\ref{thm:main_plane} is the introduction of new coordinates for the manifold $A$ and we do this in Section~\ref{Section4}. Correspondingly we give a new expression for the planar linked-twist map $\Theta$ in Section~\ref{Section5}. In Section~\ref{Section6} we introduce a new invariant cone-field and show that it is preserved by the differential $D\Theta$. Unlike the cone-field introduced by Wojtkowski, ours affords us sufficient control over the orientation of local invariant manifolds to deduce strong ergodic properties; we give the details in Section~\ref{Section7}, appealing to the work of \citeasnoun{ks} to complete the proof. We make some concluding remarks in Section~\ref{Section8}.

\section{Background to the problem}\label{Section2}
\setcounter{equation}{0}

The study of planar linked-twist maps was motivated by a number of authors. \citeasnoun{bowen} showed that certain such maps have positive topological entropy, and asked whether they possessed any ergodic properties. Similar maps were shown by \citeasnoun{braun} to arise as an approximate model of the global flow for the St\"{o}rmer problem, and were encountered by \citeasnoun{thurston} in his study of diffeomorphisms of surfaces. 

Considering briefly a more general linked-twist map $\Theta_{j,k}=\Gamma^k\circ\Phi^j$ for integers $j,k$ (where $j=k=1$ corresponds to the present case), \citeasnoun{d2} showed that if $jk\neq 0$ then there is an invariant, zero-measure Cantor set on which $\Theta_{j,k}$ is topologically conjugate to a subshift of finite type. \citeasnoun{woj} showed that under the same hypothesis, there are restrictions on the size of the annuli which guarantee that $\Theta_{j,k}$ has an ergodic partition. The restrictions are stronger for the case $jk<0$ than for the case $jk>0$; we give details for the latter case in Section~\ref{Section3}.

In an unpublished note \citeasnoun{p_preprint} considers a variety of linked-twist maps including the present kind. In \citeasnoun{p2} he shows that under certain conditions periodic saddles and homoclinic points are dense for this large class of maps and moreover that they are topologically transitive.

In recent years the study of linked-twist maps has taken on a new significance owing to developments in our understanding of the mechanisms underlying good mixing of fluids. \citeasnoun{ottbook} has shown that the single most important feature to incorporate in the design of any fluid mixing device is a `crossing of streamlines', by which we mean that flow occurs periodically in two transversal directions. That linked-twist maps provide a suitable paradigm for this design process was highlighted in \citeasnoun{ow_science} and has been discussed at much greater length in \citeasnoun{ow2} and \citeasnoun{sturman}.

This renewed emphasis on linked-twist maps in applications serves to motivate the nature of our research on this subject. While ergodic theorists have developed a very powerful and general framework for understanding the nature of ergodic behaviour in general dynamical systems (e.g.\ see the work of \citeasnoun{liverani95ergodicity}) there are very few situations relevant to applications where it is shown that the hypotheses necessary to conclude the existence of a particular ergodic property are satisfied for that particular example. This is essential for applications, and it is precisely in the spirit of our results.

The situation is very reminiscent of the development of \emph{applied} dynamical systems theory in the 1970s. Whilst it was know that generically stable and unstable manifolds of hyperbolic periodic orbits intersected transversely, and that the transverse intersections give rise to nearby Smale horseshoes, showing that this situation occurred in examples of interest to applications required significant further work (and research along these lines for concrete applications continues to this day). An excellent example illustrating this point is the work of \citeasnoun{devaney_nitecki} on showing the conditions under which the H\'{e}non map possessed an invariant set on which the dynamics was conjugate to a shift map (i.e.\ the map possessed a `horseshoe'). In that work estimates specific to the H\'{e}non map had to be carried out to show that the map satisfied the Conley-Moser criteria for the existence of such an invariant set, as given in \citeasnoun{moser}.

\section{Wojtkowski's results}\label{Section3}
\setcounter{equation}{0}

Here we describe Wojtkowski's \citeyear{woj} criteria for the planar linked-twist map to have an \emph{ergodic partition}, defined as follows:
\begin{defn}[Ergodic partition]
$\Theta$ is said to have an ergodic partition if and only if $A$ can be partitioned into at most countably many positive measure, $\Theta$-invariant, pairwise-disjoint sets $A_i$ on which the restriction of $\Theta$ is ergodic. Moreover we require that each ergodic component will be the union of finitely many Bernoulli components which are permuted by the map, i.e.\ each set $A_i$ has the form $A_i=\bigcup_{j=1}^{n(i)}A_{i,j}$ where for each $j$ the restriction of $\Theta^{n(i)}$ to $A_{i,j}$ is Bernoulli.
\end{defn}

Let $w=(u,v)\in\Sigma$ and denote by $\alpha(w)\in(0,\pi)$ the angle at which the segment connecting $w$ to $(-1,0)$ meets the segment connecting $w$ to $(1,0)$. Let
\beql
\label{eqn:alpha_condition}
\eta=\sup_{w\in\Sigma}\frac{\cot\alpha(w)}{r(w)},
\eeql
where $r(w)$ denotes the Euclidean distance from $w$ to $(-1,0)$. Recall that $c=2\pi/(r_1-r_0)$ denotes the derivative of the twist functions. Wojtkowski proved the following:
\begin{thm}[\citeasnoun{woj}] \label{thm:wojt_erg_part}
If
\beql
\label{eqn:condn_W}
c>2\eta
\eeql
then the linked-twist map $\Theta:A\to A$ has an ergodic partition.
\end{thm}
In the same paper he conjectures that under the assumptions of Theorem~\ref{thm:wojt_erg_part} then $\Theta$ also has the $K$-property. This would, by the work of \citeasnoun{ch}, imply that it has the Bernoulli property.

We discuss the proof of Theorem~\ref{thm:wojt_erg_part} briefly. It is easily argued that $\mu$-a.e.\ $w\in A$ lands in $\Sigma$ under iteration of $\Theta$ and moreover returns to $\Sigma$ infinitely many times, for those points \emph{not} satisfying this condition must be rigid rotations around one of the annuli, and must have rational angle of rotation, else their orbit would be dense and hit $\Sigma$. From the strict monotonicity of the twist function one infers that such points are contained within a set of measure zero.

Consequently for a full-measure set of points we may talk of the \emph{return map to} $\Sigma$, or just the \emph{return map} as we shall usually abbreviate it. Following Wojtkowski we define the return map on $M_+^{-1}(\Sigma)\subset L$ (rather than on $\Sigma$ itself) as the map $\Theta_{\Sigma}:M_+^{-1}(\Sigma)\to M_+^{-1}(\Sigma)$ so that
\beq*
\Theta_{\Sigma}=M_+^{-1}\circ\Theta^i\circ M_+,
\eeq*
where $i$ is the smallest (strictly) positive integer for which $\Theta^i(M_+(r,\theta))\in\Sigma$.

\begin{figure}[htp]
\centering
\includegraphics[totalheight=0.14\textheight]{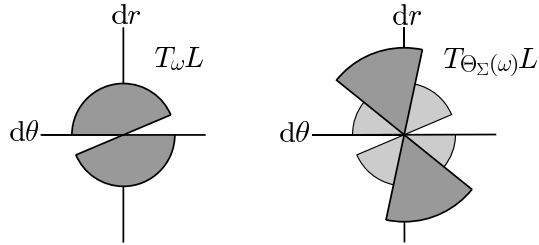}
\caption[The cone $U$]{The invariant expansive cone $U\subset T_{\omega}L$ is shown in the left-hand figure. In the right-hand figure is the image of the cone under the differential map $D\Theta_{\Sigma}$ (dark-shaded) with the original cone (light-shaded) included for comparison. Observe how the cone is mapped into itself and vectors within it are expanded.}
\label{fig:U}
\end{figure}

For $\omega=(r,\theta)\in M_+^{-1}(\Sigma)$ let $\beta_1=\text{d}r$, $\beta_2=\text{d}\theta$ give coordinates in the tangent space $T_{\omega}L$ and define the cone
\beq*
U(r,\theta)=\left\{(\beta_1,\beta_2):\frac{\beta_2}{\beta_1}\geqslant\frac{-c}{2}\right\}.
\eeq*
Wojtkowski establishes that $U$ is invariant under, and expanded by, the derivative $D\Theta_{\Sigma}$. We illustrate the situation in Figure~\ref{fig:U}. More precisely, define the cone field
\beq*
U_+=\bigcup_{(r,\theta)\in M_+^{-1}(\Sigma)}U(r,\theta)
\eeq*
and let $\|\cdot\|$ be the norm in $T_{(r,\theta)}L$ induced by the Riemannian metric, i.e.\ $\|(\beta_1,\beta_2)\|=\sqrt{\beta_1^2+r^2\beta_2^2}$. We have the following:

\begin{prop}[\citeasnoun{woj}]
\label{prop:wojtkowski_planar_cones}
$D\Theta_{\Sigma}(U_+)\subset U_+$. Furthermore there is a constant $\lambda>1$, independent of $(r,\theta)$ or $\beta$, and for vectors $\beta\in U_+$ we have $\|D\Theta_{\Sigma}\beta\|\geqslant\lambda\|\beta\|$.
\end{prop}
A detailed proof may also be found in \citeasnoun{sturman}. To arrive at the ergodic partition one determines that $\mu$-a.e.\ point returns to $\Sigma$ not just infinitely many times but with positive frequency, combines this with Proposition~\ref{prop:wojtkowski_planar_cones} to deduce non-zero Lyapunov exponents for such points and appeals to the theorem of \citeasnoun{ks}, which extends results of \citeasnoun{pes} to certain non-differentiable systems.

\section{Definition of the new coordinates}\label{Section4}
\setcounter{equation}{0}

Recall that $\mathbb{S}^1=[-\pi,\pi]$ with opposite ends identified and that we take $r_0=2$ and $r_1=\sqrt{7}$. Let $\mathcal{I}=[2,\sqrt{7}]$ and $-\mathcal{I}=[-\sqrt{7},-2]$.

\begin{figure}[htp]
\centering
\includegraphics[totalheight=0.37\textheight]{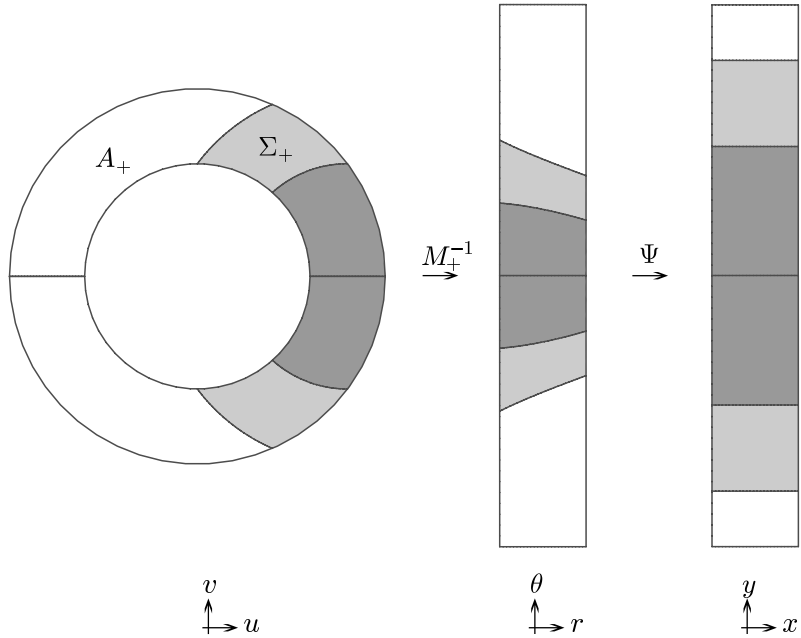}
\caption[Coordinate systems on $A+$]{The region $A_+\subset\mathbb{R}^2$, illustrated in the three coordinate systems. Left-to-right: Cartesians $(u,v)$ in the plane; polars $(r,\theta)\in\mathbb{R}_0^+\times\mathbb{S}^1$; and new coordinates $(x,y)\in\mathbb{S}^1\times\mathbb{S}^1$. Shading indicates the three regions for which $\Psi$ takes different forms, as explained in the text.}
\label{fig:new_coords}
\end{figure}

We introduce the new coordinates in two stages, starting with the annulus $A_+$. With reference to the left-hand part of Figure~\ref{fig:new_coords}, $A_+$ is divided naturally into three components: that part which intersects the annulus $A_-$ (i.e.\ the region $\Sigma$, which we have light-shaded) and the remaining connected components $A_+^{\textup{i}}$ (dark-shaded) and $A_+^{\textup{o}}$ (unshaded), which respectively lie `inside' and `outside' of the annulus $A_-$ (not shown in the figure).

We will provide first the definition and second some discussion. Let $\psi:\mathcal{I}\times[0,\pi]\to[0,\pi]$ be given by
\beq*
\psi(r,\theta)=\left\{
\begin{array}{c@{\quad}l}
{\displaystyle\frac{2\theta}{\cos^{-1}\frac{r}{4}}} &\text{if }M_+(r,\theta)\in A_+^{\textup{i}} \\ &\\
\sqrt{r^2-4r\cos\theta+4} &\text{if }M_+(r,\theta)\in\Sigma \\ &\\
\sqrt{7}+{\displaystyle\frac{\left(\pi-\sqrt{7}\right)\left(\theta-\cos^{-1}\frac{r^2-3}{4r}\right)}{\left(\pi-\cos^{-1}\frac{r^2-3}{4r}\right)}} &\text{if }M_+(r,\theta)\in A_+^{\textup{o}}
\end{array}
\right.
\eeq*
and extend $\psi$ to a function $\psi:\mathcal{I}\times\mathbb{S}^1\to\mathbb{S}^1$ by insisting that it be an odd function of $\theta$, i.e.\ that $\psi(r,-\theta)=-\psi(r,\theta)$. Finally let $\Psi:\mathcal{I}\times\mathbb{S}^1\to\mathcal{I}\times\mathbb{S}^1$ be given by
\beq*
\Psi(r,\theta)=\left(r,\psi(r,\theta)\right).
\eeq*

\begin{defn}[New coordinates $(x,y)$ on $A_+$]
Given $(u,v)\in A_+$ define $\left(x(u,v),y(u,v)\right)=\Psi\circ M_+^{-1}(u,v)$.
\end{defn}

The definitions of $\psi$ and $\Psi$ are somewhat opaque so let us now motivate them. Consider a point $(u,v)\in\Sigma_+$ and denote $(r,\theta)=M_+^{-1}(u,v)$. Then
\beq*
(x,y)=\Psi(r,\theta)=\left(r,\sqrt{r^2-4r\cos\theta+4}\right)=\left(r,\sqrt{(2-r\cos\theta)^2+r^2\sin^2\theta}\right),
\eeq*
and substituting $(u,v)=M_+(r,\theta)=(r\cos\theta-1,r\sin\theta)$ gives
\beq*
(x,y)=\left(\sqrt{(1+u)^2+v^2},\sqrt{(1-u)^2+v^2}\right).
\eeq*
So $x(u,v)$ is the Euclidean distance from $(u,v)$ to $(-1,0)$ whereas $y(u,v)$ is the Euclidean distance from $(u,v)$ to $(1,0)$; coordinates defined in this way are often called \emph{two-centre bipolar coordinates}. If we instead take $(u,v)\in\Sigma_-$ then $y$ becomes negative but still $|y|$ gives the distance to $(1,0)$.

The remainder of the definition of $\psi$ (i.e.\ for $(r,\theta)\in M_+^{-1}(A_+\backslash\Sigma)$) ensures that $\Psi:\mathcal{I}\times\mathbb{S}^1\to\mathcal{I}\times\mathbb{S}^1$ homeomorphically. This condition alone does not uniquely extend $\psi$ and so the definition given is just one of many possibilities.

We extend the new coordinates to all of $A$, making use of the rotation $\iota:\mathbb{T}^2\to\mathbb{T}^2$ given by
\beq*
\iota(x,y)=(-y,x).
\eeq*
\begin{defn}[New coordinates $(x,y)$ on $A_-$]
Given $(u,v)\in A_-\backslash\Sigma_-$ define $\left(x(u,v),y(u,v)\right)=\iota\circ\Psi\circ M_-^{-1}(u,v)$. Given $(u,v)\in\Sigma_-$ define $\left(x(u,v),y(u,v)\right)=-\iota\circ\Psi\circ M_-^{-1}(u,v)$.
\end{defn}

\begin{figure}[htp]
\centering
\includegraphics[totalheight=0.22\textheight]{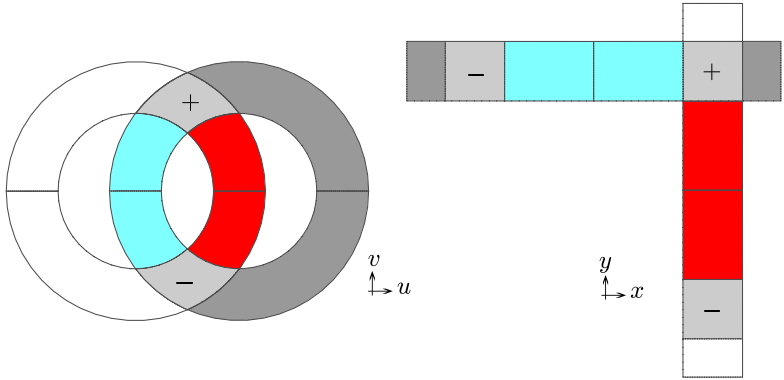}
\caption[Coordinate transformation for $A$]{The manifold $A$, illustrated in its native Cartesian coordinates and in the new coordinates $(x,y)\in\mathbb{S}^1\times\mathbb{S}^1$. Notice that there are two distinct representations of $\Sigma_-$ shown in the right-hand figure. The bottom-right representation is the correct one, the other (top-left) one is shown for sake of completion.}
\label{fig:new_coords_all}
\end{figure}

Again some discussion is required, in particular as to why we have introduced a minus sign for the case $(u,v)\in\Sigma_-$. First observe that $M_-^{-1}:A_-\to\mathcal{I}\times\mathbb{S}^1$ and $\iota:\mathbb{T}^2\to\mathbb{T}^2$ are both homeomorphisms. So \emph{without} the minus sign it follows from the definition of $\psi$ that $\iota\circ\Psi\circ M_-^{-1}:A_-\to\mathbb{S}^1\times\mathcal{I}$ homeomorphically. Moreover for  $(u,v)\in\Sigma_{\pm}$ we have
\beq*
\iota\circ\Psi\circ M_-^{-1}(u,v)=\left(\pm\sqrt{(1+u)^2+v^2},\sqrt{(1-u)^2+v^2}\right).
\eeq*
This agrees with our previous definition of $(x,y)$ on $\Sigma_+$ whereas on $\Sigma_-$ the minus sign must be introduced. Figure~\ref{fig:new_coords_all} illustrates $A$ in the original Cartesians $(u,v)$ and in the new coordinates $(x,y)\in\mathbb{S}^1\times\mathbb{S}^1$.

\section{The map expressed in the new coordinates}\label{Section5}
\setcounter{equation}{0}

We begin by showing that the linked-twist map considered (i.e.\ with $r_0=2$ and $r_1=\sqrt{7}$) has an ergodic partition; in fact we show that a whole class of linked-twist maps, including this one, have that property. The proof is simplified by the new coordinates introduced in the previous section. Following this we express $\Theta$ in the new coordinates.

\begin{lem}\label{lem:erg_part}
Let $2\leqs r_0<r_1\leqs\sqrt{7}$. Then condition (\ref{eqn:condn_W}) is satisfied.
\end{lem}
\begin{proof}
By symmetry it is enough to show that the condition holds on $\Sigma_+$. Moreover the condition is implied by
\beql
\label{eqn:new_alpha_condition}
\sup_{w\in\Sigma_+}\cot\alpha(w)<\frac{\pi}{\sqrt{7}-2}.
\eeql
Let $(x,y)=\Psi\circ M_+^{-1}(w)$ give $w\in\Sigma_+$ in the new coordinates. The angle $\alpha$ appears in a triangle in which its adjacent sides have lengths $x$ and $y$, and the opposite side has length $2$ (the Euclidean distance between the centres of the annuli). The law of cosines says that
\beql
\label{eqn:cos_alpha}
\cos\alpha=\frac{x^2+y^2-4}{2xy}.
\eeql
The partial derivative of (\ref{eqn:cos_alpha}) with respect to $x$ is given by
\beq*
\frac{\pd}{\pd x}(\cos\alpha)=\frac{1}{2y}-\frac{y^2-4}{2x^2y}
\eeq*
and, using $x,y\in\mathcal{I}$, we calculate that $1/2y\in[1/2\sqrt{7},1/4]$ and $(y^2-4)/2x^2y\in[0,3/16]$. It is easily checked that $1/2\sqrt{7}>3/16$ and so the derivative is always positive. Consequently $\cos\alpha$ is an increasing function of $x$ and so $\alpha$ is a decreasing function of $x$. By symmetry $\alpha$ is also a decreasing function of $y$. Combining these facts with (\ref{eqn:cos_alpha}) we find that $\alpha\in[\cos^{-1}(5/7),\pi/3]$. Recall that $\cot$ is positive and decreasing on $(0,\pi/2)$ so that
\beq*
\sup_{w\in\Sigma_+}\cot\alpha
=\cot\inf_{w\in\Sigma_+}\alpha
=\cot\cos^{-1}\frac{5}{7}
=\frac{5\sqrt{6}}{12}
<\frac{\pi}{\sqrt{7}-2}
\eeq*
and the proof is complete.
\end{proof}

\begin{figure}[htp]
\centering
(a)\includegraphics[totalheight=0.27\textheight]{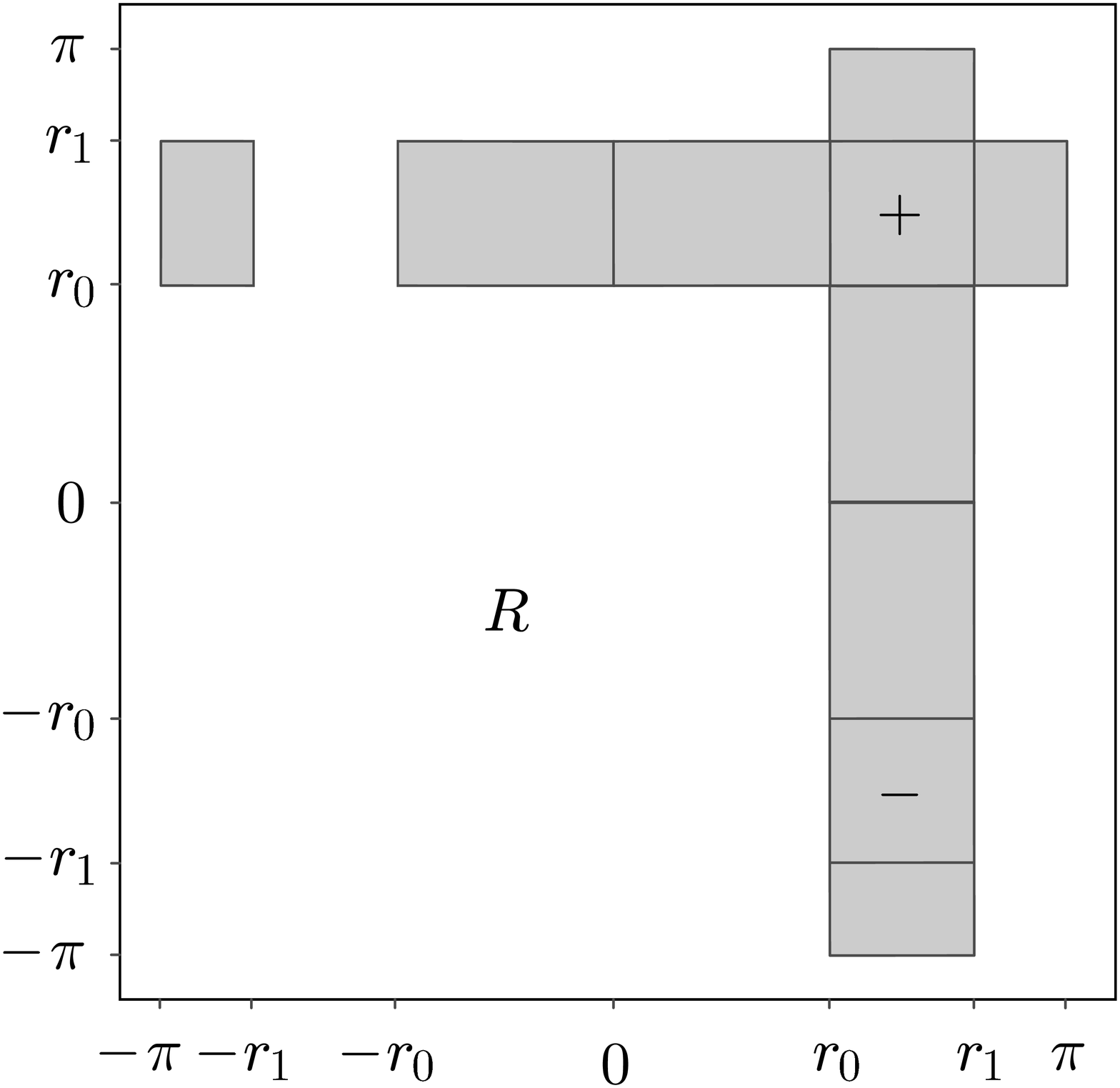}\quad
(b)\includegraphics[totalheight=0.27\textheight]{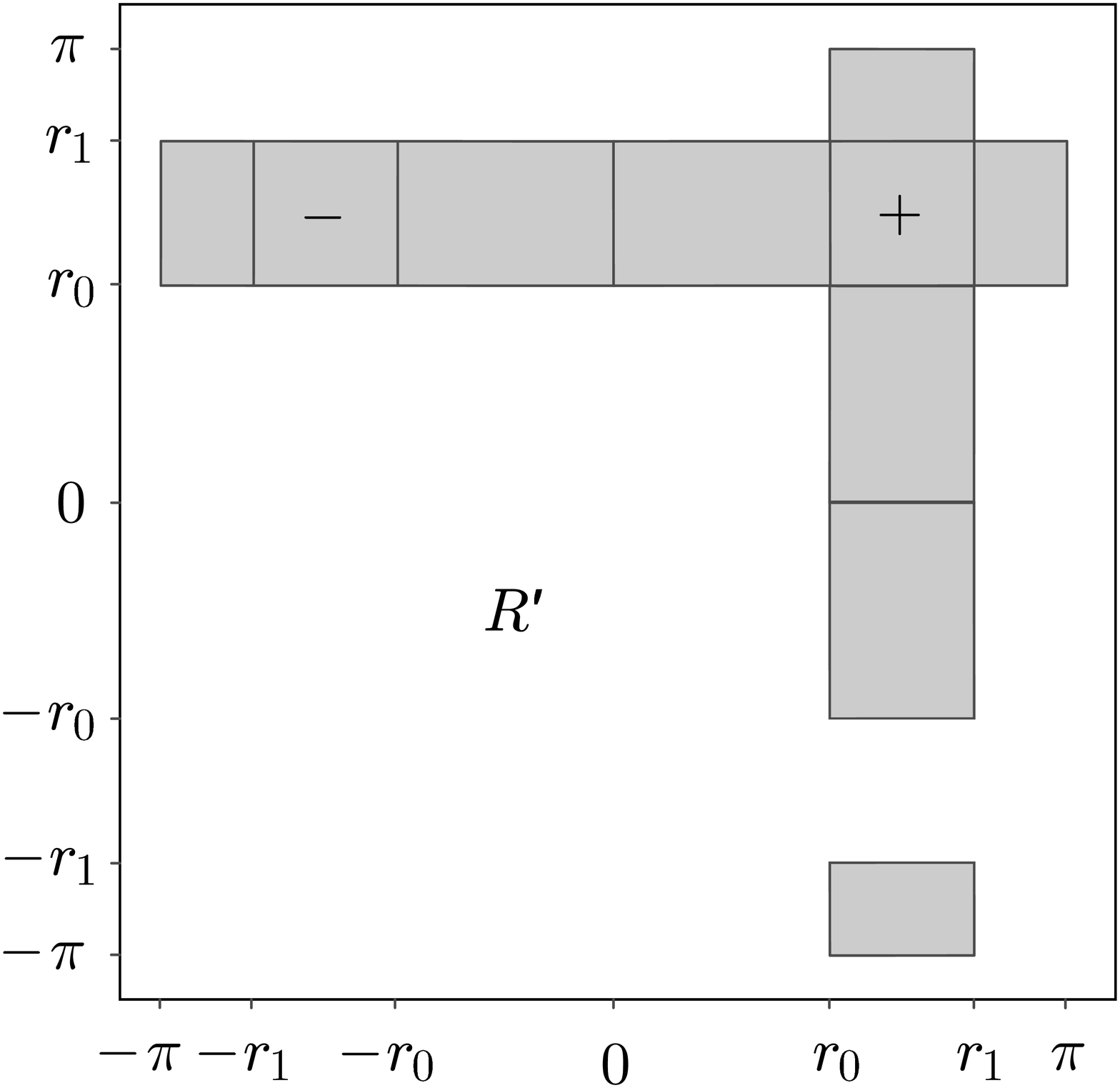}
\caption[The manifolds $R$ and $R'$]{The manifolds $R,R'\subset\mathbb{T}^2=\mathbb{S}^1\times\mathbb{S}^1$, in parts~(a) and~(b) respectively . Each is in one-to-one correspondence with $A$ but $\Sigma_-$ is represented differently in each. $F$ is a homeomorphism of $R$ and $G$ is a homeomorphism of $R'$.}
\label{fig:R_and_R'}
\end{figure}

We now express $\Theta:A\to A$ in the new coordinates. Let
\beq*
R=\left\{\Psi\circ M_+^{-1}(A_+)\right\}\cup\left\{\iota\circ\Psi\circ M_-^{-1}(A_-\backslash\Sigma_-)\right\}
\eeq*
as shown in Figure~\ref{fig:R_and_R'}(a). There is a one-to-one correspondence between points in $R$ and points in $A$. Let $F:R\to R$ denote the map $\Phi:A\to A$ in the new coordinates. Recalling our definition of $\Phi$ in (\ref{eqn:phi}) we have
\beq*\fl
F(x,y)=\Psi\circ M_+^{-1}\circ\Phi\circ M_+\circ\Psi^{-1}(x,y)=\left\{
\begin{array}{r@{\quad}l}
\Psi\circ\Lambda\circ\Psi^{-1}(x,y) &\text{if }(x,y)\in\mathcal{I}\times\mathbb{S}^1 \\
(x,y) &\text{otherwise.}
\end{array}
\right.
\eeq*
$F$ is a homeomorphism of $R$. Also let
\beq*
R'=\left\{\Psi\circ M_+^{-1}(A_+\backslash\Sigma_-)\right\}\cup\left\{\iota\circ\Psi\circ M_-^{-1}(A_-)\right\}
\eeq*
as illustrated in Figure~\ref{fig:R_and_R'}(b). Again, there is a one-to-one correspondence between points in $R'$ and points in $A$. Let $G:R'\to R'$ denote the map $\Gamma:A\to A$ in the new coordinates. Recalling our definition of $\Gamma$ in (\ref{eqn:gamma}) we have
\beq*\fl
G(x,y)=\iota\circ\Psi\circ\Lambda^{-1}\circ\Psi^{-1}\circ\iota^{-1}(x,y)=\left\{
\begin{array}{r@{\quad}l}
\iota\circ F^{-1}\circ\iota^{-1}(x,y) &\text{if }(x,y)\in\mathbb{S}^1\times\mathcal{I} \\
(x,y) &\text{otherwise.}
\end{array}
\right.
\eeq*
$G$ is a homeomorphism of $R'$.

To compose $F$ and $G$ we require the natural bijections $R\mapsto R'$, equal to $-id$ on $\mathcal{I}\times-\mathcal{I}$ and $id$ otherwise, and its inverse $R'\mapsto R$. The representation of $\Theta$ in the new coordinates is thus given by
\beq*
H=\left(R'\mapsto R\right)\circ G\circ\left(R\mapsto R'\right)\circ F=\left(R'\mapsto R\right)\circ\iota\circ F^{-1}\circ\iota^{-1}\circ\left(R'\mapsto R\right)\circ F.
\eeq*
To simplify the above expression let $\Omega,\Omega^{-1}:\mathbb{S}^1\times\mathbb{S}^1\to\mathbb{S}^1\times\mathbb{S}^1$ be given by
\beq*
\Omega^{\pm 1}(x,y)=\left\{
\begin{array}{r@{\quad}l}
\iota^{\pm 1} &\text{if }(x,y)\in\mathcal{\mathcal{I}}\times\mp\mathcal{I} \\
\iota^{\mp 1} &\text{otherwise}
\end{array}
\right.
\eeq*
then
\beq*
H=\Omega^{-1}\circ F^{-1}\circ\Omega\circ F.
\eeq*

Let $S\subset R$ be the image of the `intersection region' $\Sigma$, i.e.\ $S=\left(\mathcal{I}\times\mathcal{I}\right) \cup \left(\mathcal{I}\times-\mathcal{I}\right)$. For $z\in S$ we define the \emph{return map} $H_S:S\to S$, analogous to the return map $\Theta_{\Sigma}$.

\section{A new invariant tangent cone}\label{Section6}
\setcounter{equation}{0}

In this section we study the derivative $DH_z$ for $z\in R$. We approach this by studying the derivative $DF^{\pm 1}_z$ for $z\in R$, observing that this is the identity when $z\notin\mathcal{I}\times\mathbb{S}^1$, so the only interesting case is for $z\in\mathcal{I}\times\mathbb{S}^1$. In \emph{that} case
\beq*
F^{\pm 1}(x,y)=\Psi\circ\Lambda^{\pm 1}\circ\Psi^{-1}(x,y)=\left(x,\psi\left(x,\psi^{-1}(x,y)\pm c(x-r_0)\right)\right)
\eeq*
where $\psi^{-1}:\mathcal{I}\times\mathbb{S}^1\to\mathbb{S}^1$ is defined by $\psi^{-1}\left(x,\psi(x,y)\right)=y$. To simplify the expression write
\beq*
\tilde{y}_{\pm}=\tilde{y}_{\pm}(x,y)=\psi^{-1}(x,y)\pm c(x-r_0)\and f_{\pm}(x,y)=\psi\left(x,\tilde{y}_{\pm}\right).
\eeq*
If $D_1,D_2$ denote the usual differential operators then the Jacobians of $F^{\pm 1}$ are given by
\beq*
DF^{\pm 1}=\left( \begin{array}{cc}
1 & 0 \\ D_1f_{\pm}(x,y) & D_2f_{\pm}(x,y)
\end{array}\right).
\eeq*
The derivatives of $f_{\pm}$ are given by
\beqal
\label{eqn:D1f}
D_1f_{\pm}(x,y) &= D_1\psi(x,\tilde{y}_{\pm}) + D_2\psi(x,\tilde{y}_{\pm})\left[D_1\psi^{-1}(x,y)\pm c\right], \\
\label{eqn:D2f}
D_2f_{\pm}(x,y) &= D_2\psi(x,\tilde{y}_{\pm}) D_2\psi^{-1}(x,y).
\eeqal

\begin{figure}[htp]
\centering
\includegraphics[totalheight=0.12\textheight]{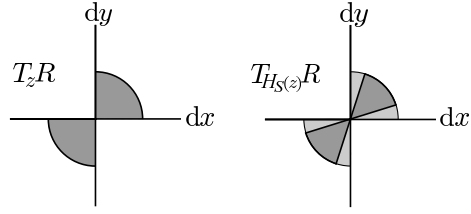}
\caption[The cone $C$]{The invariant cone $C\subset T_zR$ is shown in the left-hand figure. In the right-hand figure is the image of the cone under the differential map $DH_S$. The fact that $C$ is invariant under this differential is immediately implied by Proposition~\ref{prop:DGDF}. In contrast to the situation in Proposition~\ref{prop:wojtkowski_planar_cones} we do \emph{not} claim that this cone is expanded by $DH_S$, although it will follow from the results of Section~\ref{Section7} that this is true \emph{on average}.}
\label{fig:C}
\end{figure}

Let $b_1=\text{d}x,b_2=\text{d}y$ give coordinates in the tangent space $T_z\left(\mathbb{S}^1\times\mathbb{S}^1\right)=T_z\mathbb{T}^2$ to a point $z=(x,y)\in R$, and define the cones
\beq*
C(z)=\left\{(b_1,b_2):b_1b_2\geqs 0\right\}\and\tilde{C}(z)=\left\{(b_1,b_2):b_1b_2\leqs 0\right\}.
\eeq*
The cone $C$ is illustrated in Figure~\ref{fig:C}. Define the cone fields
\beq*
C_+=\bigcup_{z\in R}C(z)\and C_-=\bigcup_{z\in R}\tilde{C}(z).
\eeq*
The remainder of this section is devoted to proving that $DH$ preserves the cone field $C_+$.

\begin{prop}
\label{prop:DGDF}
Let $r_0=2$ and $r_1=\sqrt{7}$. If $z\in R$ and $w\in C(z)\subset T_z\mathbb{T}^2$ then
\beq*
DH_zw\in C(H(z))\subset T_{H(z)}\mathbb{T}^2.
\eeq*
\end{prop}
\begin{proof}
Some observations will simplify the task. First, from the definition of $H$ it is enough to show the result holds for each of $DF_z$ and $D\left(\Omega^{-1}\circ F^{-1}\circ\Omega\right)_z$. Second, by the chain rule and the easy observation that $D\Omega^{\pm 1}_z$ each map $C_+$ into $C_-$ and \emph{vice versa}, it is enough to show that $DF^{\pm 1}_z$ preserves $C_{\pm}$. Third, as these derivatives are otherwise the identity, it is enough to consider only $z\in\mathcal{I}\times\mathbb{S}^1$. Finally fourth, let $(b_1,b_2)\in T_z\mathbb{T}^2$ and define
\beq*
\left( \begin{array}{cc} b_1' \\ b_2'  \end{array}\right)
= DF^{\pm 1} \left( \begin{array}{cc} b_1 \\ b_2  \end{array}\right)
= \left( \begin{array}{cc} b_1 \\ b_1D_1f_{\pm}(x,y)+b_2D_2f_{\pm}(x,y) \end{array}\right).
\eeq*
We have
\beq*
\frac{b_2'}{b_1'}=D_1f_{\pm}(x,y)+\frac{b_2}{b_1}D_2f_{\pm}(x,y),
\eeq*
so it is enough to show that for each $z=(x,y)\in\mathcal{I}\times\mathbb{S}^1$ we have
\beql
\label{eqn:Df_bounds}
\pm D_1f_{\pm}(x,y)\geqs 0\quad\text{and}\quad D_2f_{\pm}(x,y)\geqs 0.
\eeql
We claim that
\beq*\fl
D_1\psi\in\left[0,\frac{7}{6}\right),\quad D_2\psi\in\left[\frac{1}{4},\sqrt{7}\right],\quad D_1\psi^{-1}\in\left[-\frac{9}{11},0\right]\and D_2\psi^{-1}\in\left[\frac{\sqrt{7}}{7},4\right].
\eeq*
It is elementary, using (\ref{eqn:D1f}) and (\ref{eqn:D2f}), to check that in this case (\ref{eqn:Df_bounds}) is satisfied. Unfortunately proving the claim will require some extensive calculations. To keep the length of the proof within reasonable limits we prove only the first assertion; the others are proved similarly. For full details see \citeasnoun{springham_thesis}. We further restrict to the case $\theta\geqs 0$; the case $\theta\leqs 0$ follows by symmetry.

Now, $D_1\psi(r,\theta)$ takes three different forms corresponding to $M_+(r,\theta)$ in each of $A_+^{\textup{i}}$, $\Sigma_+$ and $A_+^{\textup{o}}$. We deal with each case in turn.
\begin{enumerate}

\item Let $M_+(r,\theta)\in A_+^{\textup{i}}$. We will estimate the range of
\beql\label{eqn:deriv1}
\frac{\pd}{\pd r}\left(\frac{2\theta}{\cos^{-1}\frac{r}{4}}\right)=\frac{\theta}{2\left(\cos^{-1}\frac{r}{4}\right)^2\sqrt{1-\frac{r^2}{16}}}.
\eeql
The numerator attains a minimum of zero and a maximum of
\beq*
\sup_{M_+(r,\theta)\in A_+^{\textup{i}}}\theta =\sup_{r\in \mathcal{I}}\cos^{-1}\frac{r}{4} =\cos^{-1}\inf_{r\in \mathcal{I}}\frac{r}{4} =\cos^{-1}\frac{1}{2} =\frac{\pi}{3}.
\eeq*
Observe that
\beq*
\inf_{r\in \mathcal{I}}\cos^{-1}\frac{r}{4} =\cos^{-1}\sup_{r\in \mathcal{I}}\frac{r}{4} =\cos^{-1}\frac{\sqrt{7}}{4} >\frac{5}{6},
\eeq*
so the range of $\left(\cos^{-1}\frac{r}{4}\right)^2$ falls within $\left(\frac{25}{36},\frac{\pi^2}{9}\right)$. Notice also that $\frac{r^2}{16}\in\left[\frac{1}{4},\frac{1}{2}\right)$ and so $\sqrt{1-\frac{r^2}{16}}\in\left(\frac{\sqrt{2}}{2},\frac{\sqrt{3}}{2}\right]$. Combining these gives a range for (\ref{eqn:deriv1}) of $\left[0,\frac{72\pi}{150\sqrt{2}}\right)\subset\left[0,\frac{7}{6}\right)$.

\item Let $M_+(r,\theta)\in\Sigma$. We will estimate the range of
\beql\label{eqn:deriv2}
\frac{\pd}{\pd r}\left(\sqrt{r^2-4r\cos\theta+4}\right)=\frac{r-2\cos\theta}{\sqrt{r^2-4r\cos\theta+4}}.
\eeql
By design the denominator has range $\mathcal{I}$. For the numerator observe that the angle $\theta$ occurs in a triangle where the adjacent sides have lengths $2$ and $r$ and where the opposite side, call it $r'$, has length in $\mathcal{I}$. Using the law of cosines we have
\beq*
r-2\cos\theta=r-\frac{r^2+4-r'^2}{2r}=\frac{r^2-4+r'^2}{2r}.
\eeq*
The partial derivatives
\beq*
\frac{\pd}{\pd r}\left(r-2\cos\theta\right)=\frac{4+r^2-r'2}{2r^2}\and\frac{\pd}{\pd r'}\left(r-2\cos\theta\right)=\frac{r'}{r}
\eeq*
are both positive for $r,r'\in\mathcal{I}$ and so the numerator is an increasing function of each. This gives a range for the numerator of $\left[1,\frac{5\sqrt{7}}{7}\right]$ and for the quotient (\ref{eqn:deriv2}) of $\left[\frac{\sqrt{7}}{7},\frac{5\sqrt{7}}{14}\right]\subset(0,1)$.

\item Finally let $M_+(r,\theta)\in A_+^{\textup{o}}$. We estimate the range of
\beql\label{eqn:deriv3}\fl
\frac{\pd}{\pd r}\left(\sqrt{7}+\frac{\left(\pi-\sqrt{7}\right)\left(\theta-\cos^{-1}\frac{r^2-3}{4r}\right)}{\left(\pi-\cos^{-1}\frac{r^2-3}{4r}\right)}\right) =\frac{\left(r^2+3\right)(\pi-\sqrt{7})(\pi-\theta)}{r\left(\pi-\cos^{-1}\frac{r^2-3}{4r}\right)^2\sqrt{-\left(r^2-9\right)\left(r^2-1\right)}}.
\eeql
The numerator is non-negative but may be zero when $\theta=\pi$, this giving a lower bound. An upper bound requires the calculation
\beq*
\inf_{M_+(r,\theta)\in A_+^{\textup{o}}}\theta =\inf_{r\in \mathcal{I}}\cos^{-1}\frac{r^2-3}{4r} =\cos^{-1}\sup_{r\in \mathcal{I}}\frac{r^2-3}{4r} =\cos^{-1}\frac{\sqrt{7}}{7}
\eeq*
and one can check that $\pi-\theta<2$ as a consequence. Observing that $\pi-\sqrt{7}<\frac{1}{2}$ gives an upper bound for the numerator of $10$. For the denominator we need also the calculation
\beq*
\sup_{r\in \mathcal{I}}\cos^{-1}\frac{r^2-3}{4r} =\cos^{-1}\inf_{r\in \mathcal{I}}\frac{r^2-3}{4r} =\cos^{-1}\frac{1}{8}.
\eeq*
Then $\pi-\cos^{-1}\frac{r^2-3}{4r}\in\left[\pi-\cos^{-1}\frac{1}{8},\pi-\cos^{-1}\frac{\sqrt{7}}{7}\right]\subset\left(\frac{3}{2},2\right)$. Simple calculus gives $\sqrt{-\left(r^2-9\right)\left(r^2-1\right)}\in[2\sqrt{3},4]$ for $r\in\mathcal{I}$ and so the denominator takes values in $(6\sqrt{3},8\sqrt{7})$. Consequently the derivative (\ref{eqn:deriv3}) takes values in $[0,5\sqrt{3}/9)\subset[0,1)$.
\end{enumerate}
Collectively the three calculations show that $D_1\psi\in[0,\frac{7}{6})$, proving the first of the four claims. A similar approach yields the remaining parts and concludes the proof.
\end{proof}

\section{The Bernoulli property}\label{Section7}
\setcounter{equation}{0}

We now conclude the proof of our main result. For $\mu$-a.e.\ $w\in A$ \citeasnoun{sturman} show how the work of \citeasnoun{woj} gives a positive Lyapunov exponent associated to $w$. The theorem of \citeasnoun{ks} implies that there is an unstable manifold $\gamma^u(w)\subset A$ and an unstable subspace $E^u(w)\subset T_w\mathbb{R}^2$. Analogously there is a stable manifold and subspace. Moreover \citeasnoun{ks} give the following condition as sufficient for the Bernoulli property: for a.e.\ $w,w'\in A$ and for all sufficiently large natural numbers $m$ and $n$
\beql\label{eqn:rmip1}
\Theta^m\left(\gamma^u(w)\right)\cap\Theta^{-n}\left(\gamma^s(w')\right)\neq\emptyset.
\eeql

Let $z,z'$ denote $w,w'$ in the new coordinates, so $z$ is either $\Psi\circ M_+^{-1}(w)$ or $\iota\circ\Psi\circ M_-^{-1}(w)$ as appropriate and similarly for $z'$. For sake of discussion take $z=\Psi\circ M_+^{-1}(w)$ and define $\gamma^u(z)$ to be the maximal connected, smooth component of $\Psi\circ M_+^{-1}(\gamma^u(w))$ containing $z$. The piecewise smoothness of $\Psi$ ensures that there are only finitely many such smooth components and so the one containing $z$ will have positive length. By `a.e.\ $z\in R$' we mean those $z$ corresponding to some full $\mu$-measure set in $A$. We prove the following which implies (\ref{eqn:rmip1}): for a.e.\ $z,z'\in R$ and for all sufficiently large natural numbers $m$ and $n$
\beql\label{eqn:rmip2}
H^m\left(\gamma^u(z)\right)\cap H^{-n}\left(\gamma^s(z')\right)\neq\emptyset.
\eeql

There are two facts upon which our proof relies. The first is that the length, naturally defined, of $H^m\left(\gamma^u(z)\right)$ grows arbitrarily large with $m$. (Essentially this follows from Proposition~\ref{prop:wojtkowski_planar_cones} and from the one-dimensional mean-value theorem, although some care needs to be taken as $H$ is only piecewise smooth. The full proof is omitted for reasons of length, but can be found in \citeasnoun{springham_thesis}.) The second concerns the orientation of $\gamma^u(z)$:
\begin{prop}\label{prop:orient}
For a.e.\ $z\in R$ we have $E^u(z)\subset C\subset T_z\mathbb{T}^2$.
\end{prop}

\begin{proof}
Fix $z$ at which there is a positive and a negative Lyapunov exponent and consider the tangent space $T_{H^n(z)}\mathbb{T}^2$ for some $n\in\mathbb{N}$. This may be written as a direct sum $E^u_n\oplus E^s_n$ where $E^u_n=E^u(H^n(z))$ and $E^s_n=E^s(H^n(z))$ are one-dimensional unstable and stable subspaces respectively. Fix some $v_0\in C\backslash E^s_0$ then by the invariance of stable directions and of $C$ (Proposition~\ref{prop:DGDF}) one has $v_n=DH^nv_0\in C\backslash E^s_n$.

We show that if $n$ is sufficiently large then $E^u_n\subset C$ which gives the result. Assume for a contradiction that one may find arbitrarily large $n\in\mathbb{N}$ so that the inclusion does not hold. Then invariance of $C$ and of the unstable subspace imply that the it does not hold for \emph{any} $n\in\mathbb{N}$. Let
\beq*
v_n=a_ne_n^u+b_ne_n^s
\eeq*
where $e_n^{s(u)}\in E_n^{s(u)}$ is an (un)stable unit vector and where $a_n,b_n$ are non-zero and, without loss of generality, positive. Uniform expansion of unstable vectors by the return map (Proposition~\ref{prop:wojtkowski_planar_cones}) ensures that $a_n\to\infty$ as $n\to\infty$ and a similar consideration gives $b_n\to 0$. Consequently the angle between vectors $v_n$ and $e^u_n$ tends to zero in the limit, with the former strictly in $C$ and the latter by assumption not. The implication is that both approach the boundary of $C$, given by $(\d x,\d y)\in\text{span}\{(1,0)\}\cup\text{span}\{(0,1)\}$.

However $e^u_n$ cannot approach $(1,0)$, for consider some iteration $\tilde{z}=H^m(z)$ so that $H\left(\tilde{z}\right)=F\left(\tilde{z}\right)$. Then
\beq*
DH=DF=\left( \begin{array}{cc}
1 & 0 \\ D_1f_+(z) & D_2f_+(z)
\end{array}\right).
\eeq*
By continuity of the linear map $DF$ any neighbourhood $U$ of $(1,0)$ is mapped into a neighbourhood $V$ of $(1,D_1f_+(z))$. Because $D_1f_+(z)$ is bounded uniformly away from zero (as determined in the previous section), one can find $U$ so that $V\subset C$. So $e^u_n$ cannot approach $(1,0)$ as in doing so it is inevitably mapped into $C$, contradicting the assumption. Analogously one can show that $e^u_n$ cannot approach $(0,1)$ by considering some iteration $\tilde{z}=H^m(z)$ so that $H\left(\tilde{z}\right)=G\left(\tilde{z}\right)$. This gives the required contradiction.
\end{proof}

By similar arguments $E^s(z)\in\tilde{C}(z)$ and the length of $H^{-m}\left(\gamma^s(z')\right)$ diverges to infinity as $m\to\infty$. The orientations of (un)stable subspaces have the immediate consequence that gradients of (un)stable manifolds are similarly aligned on the manifold itself. The remainder of our proof is essentially geometric and is simplified by the introduction of a covering space for $R$. We do this in two stages.

Let $-R\subset\mathbb{T}^2$ denote those points $(x,y)\subset\mathbb{T}^2$ so that $(-x,-y)\in R$. Notice that $R\cap -R=\emptyset$ and let $R_1=R\cup -R$. Define $p':R_1\to R$ by $(x,y)\mapsto(x,y)$ if $(x,y)\in R$ and $(x,y)\mapsto(-x,-y)$ otherwise. Then $(R_1,p')$ is a covering space (a double cover, in fact) of $R$. The derivatives of $p'$ and its possible inverses each preserve the cones $C$ and $\tilde{C}$.

\begin{figure}[htp]
\centering
\includegraphics[totalheight=0.5\textheight]{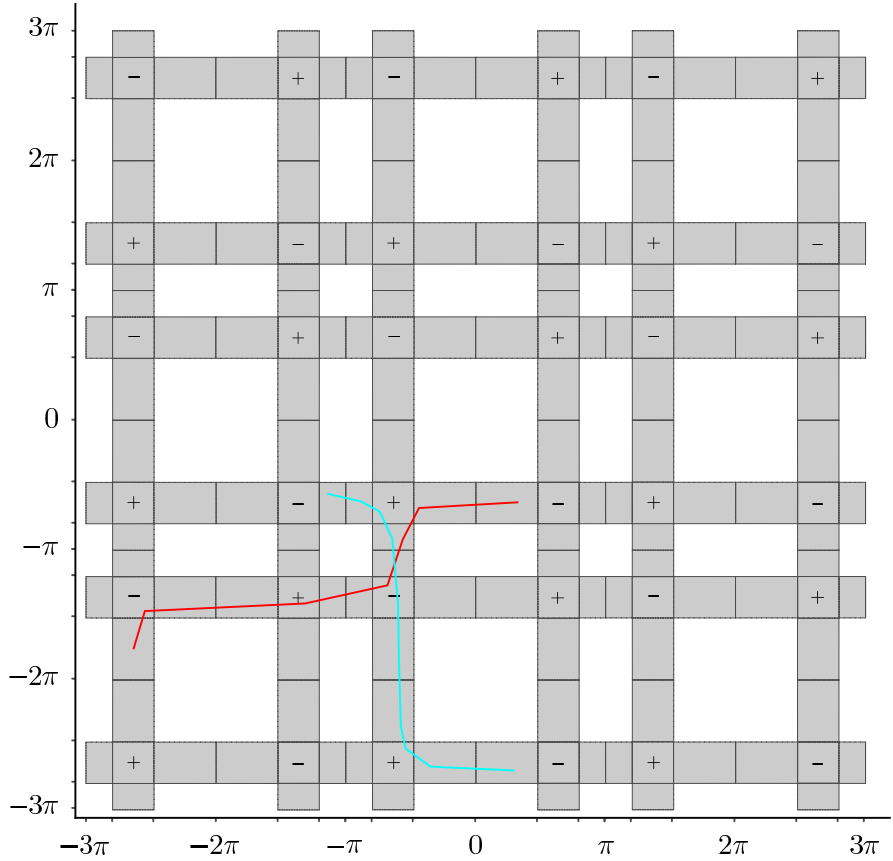}
\caption[The manifold $R_2$]{A portion of the manifold $R_2\subset\mathbb{R}^2$. Together with the map $p:R_2\to R$ this gives a covering space for $R$. In red is a typical piece of some image of a local unstable manifold for some $z\in R$. The gradient at all times is in $C$. Analogously in blue is a typical piece of some pre-image of a local stable manifold for $z'\in R$. Its gradient is in $\tilde{C}$. If each is sufficiently long then they must intersect.}
\label{fig:R_2}
\end{figure}

Let $\pi:\mathbb{R}^2\to\mathbb{T}^2$ be the natural projection which takes each coordinate modulo $2\pi$ and let
\beq*
R_2=\left\{(u,v)\in\mathbb{R}^2:\pi(u,v)\in R_1\right\}\subset\mathbb{R}^2.
\eeq*
$R_2$ has the form of a lattice, constructed by fitting together an infinite number of copies of $R_1$ and is illustrated in Figure~\ref{fig:R_2}. Let $p'':R_2\to R_1$ be the projection which takes each coordinate modulo $\mathbb{S}^1$. Then $(R_2,p'')$ gives a covering space for $R_1$. The derivatives of $p''$ and its (local) inverses preserve $C$ and $\tilde{C}$. So $(R_2,p=p'\circ p'')$ is a covering space for $R$, and $Dp$, $Dp^{-1}$ preserve $C$ and $\tilde{C}$.

It is now elementary to show that (\ref{eqn:rmip2}) is satisfied, for if we consider any sufficiently long $H^m\left(\gamma^u(z)\right)$ and any sufficiently long $H^{-n}\left(\gamma^s(z')\right)$ we can always lift them to $R_2$ in such a way that they intersect. Figure~\ref{fig:R_2} illustrates an example. The image with respect to $p$ of the intersection point is an intersection point in $R$. This completes the proof of Theorem~\ref{thm:main_plane}.

%\begin{figure}[htp]
%\centering
%\includegraphics[totalheight=0.3\textheight]{R1.eps}
%\caption[The manifold $R_1$]{The manifold $R_1=R\cup R'\subset\mathbb{T}^2$. Together with the map $p':R_1\to R$ this gives a covering space for $R$.}
%\label{fig:R_1}
%\end{figure}

\section{Concluding remarks}\label{Section8}
\setcounter{equation}{0}

Our method necessitates a strong restriction on the sizes of the annuli on which the linked-twist map is defined. The restriction is used in proving the $DH$-invariance of the tangent cone $C$ (Proposition~\ref{prop:DGDF}) where, for example, it was required to show that
\beq*\fl
D_1\psi\left(x,\psi^{-1}(x,y)+c(x-r_0)\right)+D_2\psi\left(x,\psi^{-1}(x,y)+c(x-r_0)\right)\left[D_1\psi^{-1}(x,y)+\frac{2\pi}{r_1-r_0}\right]>0,
\eeq*
for each pair $(x,y)\in [r_0,r_1]\times[0,\pi]$. Analytically determining tight estimates on the left-hand side is very difficult (even plotting it using computer algebra software requires a non-trivial effort because of the different forms taken by $\psi$ and its inverse). By comparison our approach of bounding each of $D_1\psi(\cdot,\cdot)$, $D_2\psi(\cdot,\cdot)$ and $D_1\psi^{-1}(\cdot,\cdot)$ individually is rather crude. Although the lower bound of $0$ for $D_1\psi$ is optimal, none of the other bounds established are. It is remarkable that the $DH$-invariance of $C$ may be established for \emph{any} choice of annuli using this approach and perhaps indicative that suitable bounds in fact hold for a much wider choice of annulus size. One obvious way to resolve this is to partition the domain of the functions so that tighter bounds can be established element-wise; the aforementioned computer plots might suggest a sensible partition.

We consider how far one might proceed in this manner. It is natural to wonder whether one can determine a set of values for $r_0$ and $r_1$ from which the Bernoulli property follows, and a complementary set on which it is shown not to occur. Numerical simulations of \citeasnoun{sturman} would suggest that the present method is insufficient for this task for the following reason. Essential to our ability to construct the new coordinates is that $\Sigma_{\pm}$ are disjoint, so that the shears of the respective twist maps act transversally for each point in $\Sigma$. The simulations suggest that such transversality is \emph{not} a necessary condition for good mixing. It remains an interesting open question as to whether transversality and Wojtkowski's condition (\ref{eqn:condn_W}) are \emph{sufficient} for Bernoulli.

\ack
Whilst conducting this work JS was supported by EPSRC and SRW by ONR Grant No.~N00014-01-1-0769 and EPSRC Grant EP/C515862/1. The authors are grateful to Rob Sturman for a careful reading of the text, and to Holger Waalkens and Jens Marklof for many helpful suggestions.

\section*{References}
\bibliographystyle{jphysicsB}
%\bibliography{JamesBib}

\end{document}